\documentclass[12pt]{amsart}
\usepackage{amssymb,latexsym,amsmath,amscd}

\theoremstyle{plain}
\newtheorem{thm}{Theorem}[section]
\newtheorem{prop}[thm]{Proposition}
\newtheorem{cor}[thm]{Corollary}
\newtheorem{lemma}[thm]{Lemma}

\theoremstyle{definition}
\newtheorem{defn}[thm]{Definition}
\newtheorem{defnot}[thm]{Definition/Notation}
\newtheorem{rmk}[thm]{Remark}
\newtheorem{rmks}[thm]{Remarks}
\newtheorem{notat}[thm]{Notation}

\newtheorem{cl}{Claim}
\newtheorem{ppty}{Property}
\newtheorem{ass}{Assumption}

\newcommand{\lra}{\longrightarrow}
\newcommand{\PP}{\mathbb{P}}

\newcommand{\ZZ}{\mathbb{Z}}

\newcommand{\OO}{\mathcal{O}}

\newcommand{\K}{\mathcal{K}}

\newcommand{\B}{\mathcal{B}}
\newcommand{\I}{\mathcal{I}}
\newcommand{\M}{\mathcal{M}}

\newcommand{\OP}[1]{\OO_{\PP^{#1}}}

\newcommand{\U}{\mathcal{U}}
\newcommand{\G}{\mathcal{G}}
\newcommand{\cH}{\mathcal{H}}
\newcommand{\F}{\mathcal{F}}
\newcommand{\cS}{\mathcal{S}}
\newcommand{\mmax}{{\rm max}}
\newcommand{\mmin}{{\rm min}}

\newcommand{\hgt}{\mbox{ht}\;}

\newcommand{\lt}{\mbox{LT}}
\newcommand{\lcm}{\mbox{lcm}\;}

\title{Mixed ladder determinantal varieties from two-sided ladders}

\author{Elisa Gorla}

\address{Institut f\"ur Mathematik
\\ Universit\"at Z\"urich, \hfil\break\indent Winterthurerstrasse
190, CH-8057 Z\"urich, Switzerland}

\curraddr{Max Planck Institut f\"ur Mathematik, 
\hfil\break\indent Vivatsgasse 7, D-53111 Bonn, Germany}

\email{elisa.gorla@math.unizh.ch}

\thanks{The research in this paper was carried on
while the author was a guest at the Max-Planck-Institut f\"ur 
Mathematik in Bonn. The author would like to thank
the Max Planck Foundation for its support, and for providing an 
ideal working environment.}



\begin{document}

\maketitle

{\bf Abstract:} We study the family of ideals defined by mixed size minors of 
two-sided ladders of indeterminates. 
We compute their Gr\"obner bases with respect to a skew-diagonal 
monomial order, then we use them to compute the height 
of the ideals. We show that these ideals correspond to a family of 
irreducible projective 
varieties, that we call mixed ladder determinantal varieties. 
We show that these varieties are arithmetically Cohen-Macaulay,
and we characterize the arithmetically Gorenstein ones.
Our main result consists in proving that mixed ladder determinantal 
varieties belong to the same G-biliaison class of a linear variety.

\section*{Introduction}

In this paper we study the G-biliaison class of a family of varieties,
whose defining ideals are generated by minors of a matrix of 
indeterminates. 
Other families of schemes defined by minors
have been studied in the same contest. The results obtained in this
paper are related in a very natural way to some of the results proven in~\cite{kl01},
\cite{ha04u2}, and~\cite{go05u}. In~\cite{kl01} Kleppe,
Migliore, Mir\'o-Roig, Nagel, and Peterson proved that standard
determinantal schemes are glicci, i.e. that they belong to the G-liaison class
of a complete intersection. We refer
to~\cite{kr00} for the definition of standard and good determinantal
schemes, and to~\cite{mi98b} for a comprehensive treatment of
liaison theory. Hartshorne showed in~\cite{ha04u2} 
that the double G-links produced in~\cite{kl01} can indeed
be regarded as G-biliaisons. Hence, standard determinantal schemes belong 
to the G-biliaison class of a complete intersection.
In~\cite{go05u} we defined symmetric
determinantal schemes as schemes whose saturated ideal is generated by the 
minors of size $t\times t$ of an $m\times m$ symmetric matrix with polynomial 
entries, and whose codimension is maximal for the given  $t$ and $m$.
In the same paper we proved that these schemes
belong to the G-biliaison class of a complete intersection.
In this paper we prove that mixed ladder determinantal varieties belong to the 
G-biliaison class of a linear variety, hence they are glicci.

In Section~1 we define mixed ladder determinantal ideals and varieties. In order to 
prove irreducibility and Cohen-Macaulayness of these varieties, we use some standard 
Gr\"obner bases techniques that were used e.g. in~\cite{na86}, \cite{he92a}, 
\cite{co96a}, and~\cite{go00}. 
In Theorem~\ref{gb} we compute a Gr\"obner basis of a mixed ladder determinantal 
ideal with respect to a skew-diagonal monomial order.
In Theorem~\ref{prime} we show that mixed ladder determinantal ideals are irreducible.
Hence they define irreducible varieties.
In Theorem~\ref{cm} we prove that mixed ladder determinantal varieties are 
arithmetically Cohen-Macaulay.
In Theorem~\ref{ag} we give a numerical characterization of the arithmetically Gorenstein
mixed ladder determinantal varieties.
In Theorems~\ref{codimension} and~\ref{codim} we compute the codimension of a mixed 
ladder determinantal variety by producing subsets of the ladder, whose cardinality 
agrees with the codimension of the variety. 

In Section~2 we use the results that we obtained in Section~1 in order to show that
a mixed ladder determinantal variety can be obtained from a linear variety by a 
series of ascending G-biliaisons. The explicit biliaisons are produced in the proof
of Theorem~\ref{bilmix}. They take place on suitable mixed ladder determinantal 
varieties. In Corollary~\ref{glicci} we conclude 
that mixed ladder determinantal varieties are glicci.

\section{Mixed ladder determinantal ideals and varieties}

Let $V$ be a variety in $\PP^r=\PP^r_K$, where $K$ is an algebraically
closed field. By variety we mean a reduced (not necessarily irreducible) 
scheme. We denote by $I_V$ the radical ideal
associated to $V$ in the coordinate ring of $\PP^r$. 
Let $\I_V\subset \OP{r}$ be the ideal sheaf of~$V$. Let $W$ be a 
variety containing $V$. We denote by $\I_{V|W}$ the ideal sheaf 
of $V$ restricted to $W$, i.e. the quotient sheaf $\I_V/\I_W$.
We let $H^0_*(\I)=\oplus_{s\in\ZZ}H^0(\PP^r,\I(s))$ denote \
the 0-th cohomology module of the sheaf $\I$.

In this paper we deal with ideals generated by minors of a 
matrix of indeterminates, and with the varieties that they cut out.

\begin{notat}\label{minors}
Let $X$ be an $m\times n$ matrix of indeterminates,
$$X=\left[\begin{array}{ccc}
x_{1,1} & \cdots & x_{1,n} \\
\vdots & & \vdots \\
x_{m,1} & \cdots & x_{m,n}
\end{array}\right].$$
We assume that $m\leq n$. 

Fix a choice of row indexes $1\leq i_1\leq
i_2\leq\ldots\leq i_t\leq m$ and of column indexes $1\leq j_1\leq
j_2\leq\ldots\leq j_t\leq n$. We denote by
$X_{i_1,\ldots,i_t;j_1,\ldots,j_t}$ the determinant of the $t\times t$ 
submatrix of
$X$ consisting of the rows $i_1,\ldots,i_t$ and of the columns 
$j_1,\ldots,j_t$. If all the entries involved in the minor 
$X_{i_1,\ldots,i_t;j_1,\ldots,j_t}$ belong to a subset $L$ of $X$, then we also write
$L_{i_1,\ldots,i_t;j_1,\ldots,j_t}$ for 
$X_{i_1,\ldots,i_t;j_1,\ldots,j_t}$.
\end{notat}

We are interested in the study of the ideals generated by minors that 
involve indeterminates contained in a subladder of $X$. We refer 
the reader to \cite{co93t} for information about
ladders and ladder determinantal ideals.

\begin{defnot}\label{ladder}
Let $L\subseteq X$ be a {\bf ladder}. We say that $L$ is a 
ladder if it has the following property: 
whenever $x_{i,j}, x_{k,l}\in L$ and $i\leq k, j\geq l$, then 
$x_{u,v}\in L$ for any pair $(u,v)$ with $i\leq u\leq k$ and $l\leq v\leq j$.
In other words, a ladder has the property that whenever the upper right and the 
lower left corners of a (possibly rectangular) submatrix belong to the ladder, 
then all the indeterminates in the submatrix belong to the ladder. Given 
$1=b_1< b_2<\cdots< b_h\leq m$, $1\leq a_1<\cdots<a_h=n$, 
$1\leq d_1\leq\cdots\leq d_k=m$, and $1=c_1\leq c_2\leq\cdots\leq c_k\leq n$, we let
$$L=\{x_{b,a}\;|\; b_i\leq b\leq m,\: 1\leq a\leq a_i,\: 1\leq b\leq d_j,\: 
c_j\leq a\leq n,\; \mbox{for some $i,j$}\}.$$
We call $(b_1,a_1),\ldots,(b_h,a_h)$ the {\bf upper outside corners}, and
$(d_1,c_1),\ldots,(d_k,c_k)$ the {\bf lower outside corners} of $L$. 
We allow to have more than one lower outside corner on the same row or column.
However, we assume that all the lower outside corners are distinct.
We assume that no two upper outside corners belong to the same row or to the same 
column. Since the only function of the upper outside corners is
to delimit the ladder, the assumption that no two
belong to the same row or column is not restrictive.

We include the possibility that $L=X$, i.e. {\bf a matrix is a ladder}.
We say that a ladder is {\bf one-sided} if all its lower
outside corners belong to the first column or to the last
row. Alternatively, a ladder with arbitrary lower outside corners
is one-sided if it has exactly one upper outside corner.
We regard a matrix as a one-sided ladder.
 
We let 
$$L_j=\{x_{b,a}\in L\;|\; 1\leq b\leq d_j,\: c_j\leq a\leq n\}.$$
Then $L=\cup_{j=1}^k L_j.$ 
Each $L_j$ is a one-sided ladder, with unique lower outside corner 
$(d_j,c_j)$. $L_j$ has upper outside corners 
$(b_{l_j},a_{l_j}),\ldots,(b_{h_j},a_{h_j})$, where 
$$l_j=\mmin\{i\;|\;x_{b_i,a_i}\in L_j\}\;\;\;
\mbox{and}\;\;\; h_j=\mmax\{i\;|\;x_{b_i,a_i}\in L_j\}.$$
Notice that $L_j$ is a subladder of 
$$X_j=\{x_{b,a}\in X\;|\; b_{l_j}\leq b\leq d_j, c_j\leq a\leq a_{h_j}\}.$$
Moreover, $X_j$ is the smallest submatrix of $X$ containing $L_j$.
\end{defnot}

We study ideals generated by minors of a ladder of indeterminates. 
For a given matrix $X$ and a ladder $L$, we denote 
by $I_t(L)$ the ideal generated by the $t\times t$ minors of $X$ whose entries 
belong to $L$. We call $I_t(L)$ a {\bf ladder determinantal ideal}. 
Let $K[L]$ denote the $K$-algebra 
generated by the entries of $L$. Then $I_t(L)\subseteq K[L]$. Narasimhan showed 
in~\cite{na86} that these ideals are prime, hence they define irreducible projective 
varieties.

\begin{defn}
Let $X$ be a matrix of indeterminates, let $L\subseteq X$ be a ladder. 
Fix a positive integer $t$.
We call {\bf ladder determinantal variety} a variety with ideal $I_t(L)$. 
We regard it as a subvariety of $\PP^r=Proj(K[L])$.
\end{defn}

As in~\cite{go00}, we extend the definition of ladder determinantal ideal 
to allow minors of different 
size in different regions of the ladder. In the paper of Gonciulea and Miller, 
they call mixed ladder determinantal varieties those corresponding to mixed 
ladder determinantal ideals coming from one-sided ladders. Moreover, they assume 
that no two consecutive outside corners belong to the same
row or column (notice that their upper corners correspond to our
lower corners, and viceversa). 
In this paper we drop both of those assumptions, so we deal with a larger 
family of varieties than in~\cite{go00}.

\begin{defn}
Let $X$ be a matrix of indeterminates, let $L\subseteq X$ be a ladder as in 
Definition~\ref{ladder}.
Fix a vector of positive integers $t=(t_1,\ldots,t_k)$. A {\bf mixed ladder 
determinantal ideal} is an ideal of the form 
$I_t(L):=I_{t_1}(L_1)+\cdots+I_{t_k}(L_k).$
We call {\bf mixed ladder determinantal variety} a variety defined by a mixed 
ladder determinantal ideal. As for ladder determinantal varieties, if 
$I_V=I_t(L)$ then $V$ is a subvariety of $\PP^r=Proj(K[L])$.
\end{defn}

\begin{rmk}\label{nndg}
In Theorem~\ref{prime} we show that every mixed ladder determinantal ideal is prime.
Therefore every mixed ladder determinantal ideal is the ideal of a mixed ladder 
determinantal variety. Moreover, mixed ladder determinantal varieties are 
irreducible. 
\end{rmk}

\begin{rmk}
Every ladder determinantal variety is a mixed ladder 
determinantal variety, whose vector $t$ has all equal entries.
\end{rmk}

We adopt the following conventions on the ladder $L$. 
All of the assumptions can be made without loss of generality.

\begin{ass}\label{nondeg}
All the entries of $L$ are involved in some $t$-minor. Otherwise we delete 
the superfluous entries. It is easy to check that the resulting set is again a 
ladder.
\end{ass}
\begin{ass}\label{corners}
No two consecutive corners coincide; no two upper outside corners belong to the
same row or to the same column.
\end{ass}
\begin{ass}\label{assumpt}
$I_{t_i}(L_i)\not\subseteq I_{t_j}(L_j)$ for $i\neq j$. Therefore
$$d_{j+1}-d_j>t_{j+1}-t_j\;\;\;\mbox{and}\;\;\; c_{j+1}-c_j>t_j-t_{j+1}.$$
In particular, if $c_j=c_{j+1}$ then $t_j<t_{j+1}$. If $d_j=d_{j+1}$, then 
$t_j>t_{j+1}$.
\end{ass}

We allow disconnected ladders, i.e. ladders for which $L_i\cap L_{i+1}=\emptyset$ for 
some $i$. We also allow ladders that are contained in a proper submatrix of $X$ (i.e. 
we do not distinguish between mixed ladder determinantal varieties and cones over them). 
This is necessary in our context, since in some of the arguments we will need to delete 
some entries of a ladder (see e.g. Proposition~\ref{maxmin} and Theorem~\ref{bilmix}). 
So even if we start with a ladder $L$ that is connected and that is not contained in any
proper submatrix of $X$, these properties will not always be preserved.

Notice that Assumption~\ref{nondeg} is not necessarily verified by the subladders $L_j$, 
even if it is verified by $L$. Assumption~\ref{corners} is automatically verified by the 
subladders $L_j$. Assumption~\ref{assumpt} is trivially verified by every $L_j$, since 
they have exactly one lower outside corner $(d_j,c_j)$ and all the minors have the same 
size $t_j$.

In Section~2 we show that mixed ladder determinantal varieties belong to 
the G-biliaison class of a complete intersection.
Prior to that, we need to establish some properties of mixed ladder determinantal 
varieties, such as irreducibility and Cohen-Macaulayness. Irreducibility of ladder 
determinantal varieties was proven by Narasimhan in \cite{na86}, while 
irreducibility of mixed ladder determinantal varieties associated to one-sided 
ladders was proven by Gonciulea and Miller in \cite{go00}. Cohen-Macaulayness of 
ladder determinantal varieties was established by Herzog and Trung in \cite{he92a}.
The analogous result  for mixed ladder determinantal varieties associated to 
one-sided ladders was proven by Gonciulea and Miller in \cite{go00}.

In this section we compute the codimension of mixed ladder determinantal 
varieties. Analogous results were proven by Herzog and Trung (for ladder determinantal 
varieties) and by Gonciulea and Miller (for mixed ladder determinantal varieties 
associated to one-sided ladders) in the papers mentioned above. In order to compute 
the codimension, 
we compute a Gr\"obner basis of mixed ladder determinantal ideals with respect to a 
skew-diagonal monomial order.

\begin{notat}
We consider the lexicographic order on the monomials of $R$, induced by the total
order $<$ on the variables: 
$$x_{d,c}<x_{b,a}\;\;\;\mbox{if}\;\; b<d\;\;\;\mbox{or}\;\; b=d,\; a>c.$$
The leading term of a minor with respect to this term order is the product of the 
elements on its skew-diagonal.

For $f\in K[L]$ we denote by $\lt(f)$ the leading term of $f$ with 
respect to $<$. For an ideal $I\subset K[L]$ we denote $\lt(I)$ the leading term ideal
of $f$ with respect to $<$.
\end{notat}

We now compute a Gr\"obner basis of $I_t(L)$ with respect to the monomial 
order $<$. We will use this Gr\"obner basis to compute the height
of the mixed ladder determinantal ideal $I_t(L)$.
The argument is a standard one, see for example \cite{co96a} and \cite{go00}.
We first recall two known results.

\begin{lemma}\label{gbint}(\:\cite{co93t}, Lemma~1.3.13)
Let $I\subset K[X]$ be an ideal with Gr\"obner basis $G$. Let $L\subseteq X$
have the property that if $f\in G$ and $\lt(f)\in K[L]$, then $f\in K[L]$.
Then $G\cap K[L]$ is a Gr\"obner basis of $I\cap K[L]$.
\end{lemma}

\begin{lemma}\label{gblemma}(\:\cite{co93t}, Lemma~1.3.14)
Let $I,J\subset K[x_1,\ldots,x_r]$ be homogeneous ideals. Let $F$ be a 
Gr\"obner basis of $I$ and $G$ be a Gr\"obner basis of $J$. Then: $F\cup G$ 
is a Gr\"obner basis of $I+J$ if and only if for all $f\in F$ and $g\in G$ there is an 
$h\in I\cap J$ such that $\lt(h)=\lcm(\lt(f),\lt(g))$.
\end{lemma}

In the next theorem we produce a Gr\"obner basis of $I_t(L)$ with respect to $<$.

\begin{thm}\label{gb}
Let $\G_j$ be the set of determinants of the $t_j\times t_j$ submatrices of 
$L_j$, such that at most $t_i-1$ rows belong to $L_i$ for $1\leq i<j$ and 
at most $t_i-1$ columns belong to $L_i$ for $j\leq i<k.$   
Let $\G=\cup_{j=1}^k \G_j$. Then
$\G$ is a Gr\"obner basis of $I_t(L)$ with respect to $<$.
\end{thm}

\begin{proof}
The $t_j\times t_j$ minors of $L_j$ that do not belong to $\G_j$ 
can be written as a combination of the minors of $\G_i$ for some $i\neq j$. 
Therefore, $\G$ generates $I_t(L)$.
Let $\F_j$ be the set of $t_j\times t_j$ minors in $L_j$, $\F_j\supseteq\G_j$.
Moreover for each $f\in\F_j\setminus\G_j$, $\lt(f)$ is a multiple of $\lt(g)$ 
for some $g\in\G_i$, $i\neq j$. Let $\F=\cup_{j=1}^k\F_j$. 
From what we just observed, if $\F$ is a Gr\"obner basis of $I_t(L)$, then $\G$
is a Gr\"obner basis of $I_t(L)$ as well.

In order to prove that $\F$ is a Gr\"obner basis, 
we proceed by induction on $k$. If $k=1$ then $L$ is a one-sided ladder, $\F=\G$,
and the thesis follows from Corollary~3.3 of~\cite{na86}. 

If $k=2$, then $L=L_1\cup L_2$.
Let $(d_1,1),(m,c_2)$ be the lower outside corners of $L$.
Let $(b_1,1),(b_2,a_2),\ldots,(b_l,a_l)$ be the upper outside corners of $L_1$ and 
let $(b_g,a_g),\ldots,$ $(b_{h-1},a_{h-1}),(b_h,n)$ be the upper outside corners 
of $L_2$.
Let $X_1$ and $X_2$ be as in Notation~\ref{ladder}.
Let $\cH_i$ be the set of $t_i\times t_i$ minors contained in $X_i$, $i=1,2$.
By Theorem~4.5.4 in~\cite{go00}, $\cH_1\cup\cH_2$ is a Gr\"obner basis of 
$I=I_{t_1}(X_1)+I_{t_2}(X_2)$.
By Lemma~\ref{gbint} $(\cH_1\cup\cH_2)\cap K[L]$ is a Gr\"obner 
basis of $I\cap K[L]$. The assumption of the lemma holds because 
the leading term of each monomial is the product of the elements on the 
skew-diagonal, and since $L$ is a ladder. We claim that 
$\F_1\cup\F_2=(\cH_1\cup\cH_2)\cap K[L]$. In fact, 
$$\F_i=\cH_i\cap K[L_i]=\cH_i\cap K[L]$$
for $i=1,2$. Since $\F=\F_1\cup\F_2$ generates $I_t(L)$, it is a 
Gr\"obner basis of it.

Assume that the thesis is true for $k-1$. Let $L'=\cup_{j=1}^{k-1}L_j$; $L'$ has
lower outside corners $(d_1,c_1),\ldots, (d_{k-1},c_{k-1})$. Let 
$t'=(t_1,\ldots,t_{k-1})$. Then $$L=L'\cup L_k\;\;\;\mbox{and}\;\;\; 
I_t(L)=I_{t'}(L')+I_{t_k}(L_k).$$ $\F_k$ is a Gr\"obner basis of 
$I_{t_k}(L_k)$ by Corollary~3.3 of~\cite{na86}. 
Let $\F'=\cup_{j=1}^{k-1}\F_j$. By induction hypothesis $\F'$
is a Gr\"obner basis of $I_{t'}(L')$. Moreover, $\F=\F'\cup\F_k$. 
Let $f\in\F'$ and $g\in\F_k$. Then $f\in\F_j$ for some $1\leq j\leq k-1$.
Let $H=L_j\cup L_k$. We have seen above that $\F_j\cup\F_k$ is a Gr\"obner basis 
of $I_{(t_j,t_k)}(H)=I_{t_j}(L_j)+I_{t_k}(L_k)$. Then by Lemma~\ref{gblemma}
there is an $h\in I_{t_j}(L_j)\cap I_{t_k}(L_k)\subseteq I_{t'}(L')
\cap I_{t_k}(L_k)$ such that $\lt(h)=\lcm(\lt(f),\lt(g))$. 
Therefore, $\F=\F'\cup\F_k$ is a Gr\"obner basis of $I_t(L)$, again 
by Lemma~\ref{gblemma}.
\end{proof}

We now use Theorem~\ref{gb} to compute the height of 
mixed ladder determinantal ideals. In Theorem~\ref{codim} we produce a subladder $L'$ 
of $L$ such that the number of entries in $L'$ equals the height of $I_t(L)$.
We start by fixing the notation.

\begin{notat}
For $1\leq i\leq j\leq k$ we define 
$$L_{i,j}=(L_i\setminus L_{i-1})\cap (L_j\setminus L_{j+1}).$$
We adopt the convention that $L_{-1}=L_{k+1}=\emptyset$.
Each $L_{i,j}$ is a one-sided ladder with exactly one lower outside corner. 
Some of the $L_{i,j}$'s may be empty.  
Set $$t_{i,j}=\mmin\{t_l\;|\; L_{i,j}\subseteq L_l\}.$$
If $L_{i,j}=\emptyset$, we set $t_{i,j}=1$.
\end{notat}

\begin{rmk}
It is easy to check that $\cup_{i,j}L_{i,j}=L$, and that
$L_{i,j}\cap L_{p,q}\neq\emptyset$ if and only if $(i,j)=(p,q)$. 
\end{rmk}

Moreover, one has the following inequalities between the $t_{i,j}$'s.

\begin{lemma} With the notation above, one has
$$t_{i,j-1}\geq t_{i,j}\;\;\;\mbox{and}\;\;\; t_{i,j}\geq t_{i-1,j}$$
whenever $L_{i,j},L_{i,j-1},L_{i-1,j}\neq\emptyset.$ For all $i,j$ one has 
$t_{i,j}\geq 1.$
\end{lemma}

\begin{proof}
$t_{i,j}\geq 1$ for all $i,j$ since $t_l\geq 1$ for $l=1,\ldots,k$. 
If $L_{i,j}=\emptyset$, then $t_{i,j}=1$ by definition.
The other inequalities follow from the observation that 
$$L_{i,j}=(L_i\setminus\bigcup_{l=1}^{i-1}L_l)\bigcap
(L_j\setminus\bigcup_{l=j+1}^k L_l).$$ 
\end{proof}

We are now ready to construct a subset $\B\subseteq L$, whose cardinality 
agrees with the dimension of the mixed ladder determinantal variety defined by 
$I_t(L)\subset K[L]$.

\begin{thm}\label{codimension}
Let $\B_{i,j}$ be the subset of $L_{i,j}$ consisting of the elements that belong 
to the first $t_{i,j-1}-t_{i,j}$ columns or to the last $t_{i,j}-t_{i-1,j}$ rows. 
Let $\B=\cup_{i,j}\B_{i,j}$. Then $$\hgt I_t(L)=|L\setminus\B|.$$
\end{thm}

\begin{proof}
It is well known that $\hgt I_t(L)=\hgt \lt(I_t(L)).$ Moreover 
$\hgt \lt(I_t(L))=|\cS|$, where $\cS$ is a subset of the indeterminates in $L$ 
of minimal cardinality among those with the 
\begin{ppty}\label{S} 
Each monomial in a 
set of generators of $\lt(I_t(L))$ contains a variable from $\cS$. 
\end{ppty}
Geometrically, $\lt(I_t(L))$ corresponds to an intersection $\U$ of linear spaces, 
and the indeterminates in $\cS$ cut out a linear space of minimal codimension 
among those which contain $\U$ (that therefore has the same codimension as $\U$).
We claim that $L\setminus \B$ is a subset of $L$ with the Property~\ref{S}. In fact, 
we can choose as a set of generators for $\lt(I_t(L))$ the set of all 
skew-diagonal terms of $t_j\times t_j$ minor in $L_j$. Notice moreover that 
$t_j\geq t_{i,j}$ and that $L_{i,j}\setminus\B_{i,j}\neq\emptyset$ if 
$L_{i,j}\neq\emptyset$.

In order to complete the proof, we need to show that any  $\cS'$ with 
$|\cS'|<|L\setminus\B|$ does not satisfy the Property~\ref{S}. Let 
$\cS=L\setminus\B$, and let
$$\delta_u=\{x_{b,a}\in L\;|\; b+a=u+1\},\;\; u\geq 1.$$
Then $L=\cup_{u\geq 1}\delta_u$ and each two distinct elements in the union
are disjoint. Since 
$$\sum_{u\geq 1}|\cS\cap\delta_u|=|\cS|>|\cS'|=\sum_{u\geq 1}|\cS'\cap\delta_u|$$
then $|\cS\cap\delta_u|>|\cS'\cap\delta_u|$ for some $u\geq 1.$ 
Notice that if $\delta_u\cap L_j\neq\emptyset$, then $\delta_u\subseteq L_j$.
This follows from the definition of ladder. So let $j$
be maximum with the property that $\cS\cap\delta_u\subseteq L_j$. Then
$$|\cS'\cap\delta_u\cap L_j|\leq |\cS'\cap\delta_u|<|\cS\cap\delta_u|=
|L_j\cap\delta_u|-|(\cS\setminus L_j)\cap\delta_u|=|L_j\cap\delta_u|-(t_j-1).$$
Therefore there are $t_j$ distinct variables in 
$(\delta_u\cap L_j)\setminus\cS'$. Their product corresponds to a minimal 
generator of $\lt(I_t(L))$ that does not involve any element of $\cS'$.
\end{proof}

We now construct a subladder $L'$ of $L$ such that the number of entries in $L'$ 
equals the height of $I_t(L)$.

\begin{thm}\label{codim}
Let $L$ be a ladder with upper outside corners and lower outside corners as
in Notation~\ref{ladder}. Let $L'$ be the sublattice with the same upper 
outside corners as $L$, and lower outside corners 
$\{(d_j-t_j+1,c_j+t_j-1)\;|\; 1\leq j\leq k\}$. Then $$\hgt I_t(L)=|L'|.$$ 
\end{thm}

\begin{proof}
First notice that $L'$ is a ladder, since the coordinates of its outside 
lower corners verify the inequalities 
\begin{equation}\label{dct}
d_{j+1}-t_{j+1}+1>d_j-t_j+1\;\;\;\mbox{and}\;\;\; c_{j+1}+t_{j+1}-1>c_j+t_j-1
\end{equation}
by Assumption~\ref{assumpt}. Notice that it follows from~(\ref{dct}) that no two 
outside lower corners of $L'$ belong to the same row or column.
Let $\M$ be the set of minors in $\G$ that involve consecutive rows and columns
(see the statement of Theorem~\ref{gb} for the definition of $\G$). 
\begin{cl}
$|\M|=|L\setminus\B|$.
\end{cl}
There is a bijection between the two sets, given by the fact that every minor in 
$\M$ corresponds uniquely to its lower left corner.
\begin{cl}
$|\M|=|L'|$.
\end{cl} 
Again, there is a bijection between the two sets. Every minor in $\M$ corresponds 
uniquely to its upper right corner. The two claims together with 
Theorem~\ref{codimension} prove the thesis. 
\end{proof}

\begin{rmk}
Theorem~\ref{codim} also follows from the observation that for each 
$1\leq i\leq h$ and 
$1\leq j\leq k$ such that $L_{i,j}\neq\emptyset$, we delete
$r_{i,j}=t_{i,j}-t_{i-1,j}$ rows and $c_{i,j}=t_{i,j-1}-t_{i,j}$
columns from $L_{i,j}$. Therefore, the total number of entries deleted
from a row that belongs to $L_i\setminus L_{i-1}$ is 
$$\sum_{j: \: L_{i,j}\neq\emptyset} c_{i,j}=\sum_{j: \: L_{i,j}\neq\emptyset} 
(t_{i,j}-t_{i-1,j})=t_i-1.$$
Similarly, the total number of entries deleted
from a column that belongs to $L_j\setminus L_{j+1}$ is 
$$\sum_{i: \: L_{i,j}\neq\emptyset} r_{i,j}=\sum_{i: \:
  L_{i,j}\neq\emptyset} (t_{i,j-1}-t_{i,j})=t_j-1.$$
\end{rmk}

\begin{rmk}
In Theorem~\ref{bilmix} we show that if $V$ is a mixed ladder determinantal 
variety with $I_V=I_t(L)$, then $V$ is G-bilinked in $t_1+\cdots+t_k-k$ steps 
to the linear variety associated to $I_1(L')$.
\end{rmk}

Following the proof of Theorem~4.1 in~\cite{na86}, one can easily show that mixed ladder 
determinantal ideals are prime. Hence they define irreducible projective varieties.

\begin{thm}\label{prime}
Let $X_i\subseteq X$ be as in Notation~\ref{ladder}. 
Let $X'=X_1\cup\ldots\cup X_k$ and let 
$I_t(X')=I_{t_1}(X_1)+\cdots+I_{t_k}(X_k)\subseteq K[X']$. 
Then $$I_t(X')\cap K[L]=I_t(L).$$
In particular, $I_t(L)$ is a prime ideal. 
\end{thm}

\begin{proof}
It is enough to show that $I_t(X')\cap K[L]\subseteq I_t(L).$
For $j=1,\ldots,k$, let $\G_j$ be the set of determinants of the 
$t_j\times t_j$ submatrices of $X_j$, such that at most $t_i-1$ rows belong to $X_i$ for 
$1\leq i<j$ and at most $t_i-1$ columns belong to $X_i$ for $j\leq i<k.$   
Let $\G=\cup_{j=1}^k \G_j$. By Theorem~\ref{gb}
$\G$ is a Gr\"obner basis of $I_t(X')$ with respect to $<$.
Let $f\in I_t(X')\cap K[L]$ be a homogeneous form. Then 
$$f=\lt(f)+g=\lt(m)h+g$$
where $g=f-\lt(f)$, $\lt(m)$ is the leading term of some element $m\in\G$, and $h\in K[L]$ 
is a monomial (possibly multiplied by a constant). Then $$f=mh+f'$$
where $f'=f-mh\in I_t(X')\cap K[L]$, $deg(f)=deg(f')$ and $\lt(f')<\lt(f)$.
The proof is complete by induction on $\lt(f)$. 
\end{proof}

Any mixed ladder determinantal variety has a nonempty open subset which is isomorphic 
to a nonempty open subset of a ladder determinantal variety. 
We prove this in Proposition~\ref{maxmin}, by repeatedly localizing the coordinate ring of 
the variety at a corner. We need the following lemma.

\begin{lemma}\label{local}
Let $L$ be a ladder and let $I_t(L)$ be a mixed ladder determinantal ideal,  
as in Notation~\ref{ladder}. Let $d_0=0$ and $c_{k+1}=n+1$.
Fix a $j\in\{1,\ldots,k\}$ such that
$$t_j>1,\;\;\; d_{j-1}<d_j,\;\;\; c_{j+1}>c_j.$$
Let $M$ be the ladder obtained from $L$ by removing the indeterminates 
$$x_{d_{j-1}+1,c_j},\ldots,x_{d_j,c_j},x_{d_j,c_j+1},\ldots,x_{d_j,c_{j+1}-1}.$$
The ladder $M$ has lower outside corners 
$$\{(d_1,1),\ldots,(d_{j-1},c_{j-1}),(d_j-1,c_j+1),(d_{j+1},c_{j+1}),\ldots,
(m,c_k)\}$$ and the same upper outside corners as $L$.
Let $t'=(t_1,\ldots,t_{j-1},t_j-1,t_{j+1},\ldots,t_k).$
Then there is an isomorphism $$K[L]/I_t(L)[x_{d_j,c_j}^{-1}]\cong 
K[M]/I_{t'}(M)[x_{d_{j-1}+1,c_j},\ldots,x_{d_j,c_j},x_{d_j,c_j+1},\ldots,x_{d_j,c_{j-1}-1}]
[x_{d_j,c_j}^{-1}].$$
\end{lemma}

\begin{proof}
Define $$\varphi:K[L][x_{d_j,c_j}^{-1}]\lra 
K[M][x_{d_{j-1}+1,c_j},\ldots,x_{d_j,c_j},x_{d_j,c_j+1},\ldots,x_{d_j,c_{j+1}-1}]
[x_{d_j,c_j}^{-1}]$$ by the assignment 
$$\varphi(x_{u,v})=\left\{\begin{array}{ll}
x_{u,v}+x_{u,c_j}x_{d_j,v}x_{d_j,c_j}^{-1} & 
\mbox{if $u\neq d_j, v\neq c_j$ and $x_{u,v}\in L_j$,}\\
x_{u,v} & \mbox{otherwise.}
\end{array}\right.$$
The inverse of $\varphi$ is $$\psi:K[M]
[x_{d_{j-1}+1,c_j},\ldots,x_{d_j,c_j},x_{d_j,c_j+1},\ldots,x_{d_j,c_{j+1}-1}]
[x_{d_j,c_j}^{-1}]\lra K[L][x_{d_j,c_j}^{-1}]$$
defined by the assignment 
$$\psi(x_{u,v})=\left\{\begin{array}{ll}
x_{u,v}-x_{u,c_j}x_{d_j,v}x_{d_j,c_j}^{-1} & 
\mbox{if $u\neq d_j, v\neq c_j$ and $x_{u,v}\in M_j$,}\\
x_{u,v} & \mbox{otherwise.}
\end{array}\right.$$
Let $$J_t(L)=I_t(L)K[L][x_{d_j,c_j}^{-1}]$$ and 
$$J_{t'}(M)=I_{t'}(M)K[M]
[x_{d_{j-1}+1,c_j},\ldots,x_{d_j,c_j},x_{d_j,c_j+1},\ldots,x_{d_j,c_{j+1}-1}]
[x_{d_j,c_j}^{-1}].$$
It suffices to show that $\varphi(J_t(L_j))=J_{t'}(M_j)$ for $1\leq j\leq k$, 
where $M_j=L_j\cap M$. Since $L_j$ is a one-sided ladder, the thesis follows 
from Lemma~5.1.4 in~\cite{go00}.
\end{proof}

\begin{prop}\label{maxmin}
Let $I_t(L)$ be a mixed ladder determinantal ideal as in Notation~\ref{ladder}.
Let 
$$t_{\mmax}=\mmax\{t_1,\ldots,t_k\},\;\;\;\;\; t_{\mmin}=\mmin\{t_1,\ldots,t_k\}.$$
Let $L_{\mmax}$ be the ladder obtained from $L$ by enlarging each $L_j$ along its lower 
border by a strip of $t_{\mmax}-t_j$ new indeterminates. In other words, $L_{\mmax}$ 
has the same upper outside corners as $L$, and lower outside corners
$\{(d_j+(t_{\mmax}-t_j),c_j-(t_{\mmax}-t_j))\}.$
Let $L_{\mmin}$ be the ladder obtained from $L$ by removing from each $L_j$ along its lower 
border a strip of $t_j-t_{\mmin}$ indeterminates. In other words, $L_{\mmin}$ 
has the same upper outside corners as $L$, and lower outside corners
$\{(d_j-(t_j-t_{\mmin}),c_j+(t_j-t_{\mmin}))\}.$
Let $I_{\mmax}=I_{t_{\mmax}}(L_{\mmax})$ and $I_{\mmin}=I_{t_{\mmin}}(L_{\mmin})$.
Then there are $u\in K[L_{\mmax}]/I_{t_{\mmax}}(L_{\mmax})$, 
$v\in K[L]/I_t(L)$, and indeterminates $w_1,\ldots,w_i,z_1,\ldots,z_j$ such that
$$K[L_{\mmax}]/I_{\mmax}[u^{-1}]\cong K[L]/I_t(L)[w_1,\ldots,w_i]
[w_1^{-1},\ldots,w_{i_0}^{-1}]$$
and 
$$K[L_{\mmin}]/I_{\mmin}[z_1,\ldots,z_j][z_1^{-1},\ldots,z_{j_0}^{-1}]\cong 
K[L]/I_t(L)[v^{-1}]$$
for some $i_0\leq i$ and $j_0\leq j$.
\end{prop}

\begin{proof}
The thesis easily follows from Lemma~\ref{local}. We sketch how it can be derived 
from it.
In order to prove the first of the two ring isomorphisms,
keep localizing the ring $K[L_{\mmax}]/I_{\mmax}$ at the 
indeterminates in the appropriate lower outside corners of $L_{\mmax}$. 
The size of the minors considered in each region decreases by 1 at each 
localization, and by Lemma~\ref{local} we have an isomorphism as claimed. 
A similar argument can be made starting 
from $K[L]/I_t(L)$. The products of the elements that were inverted give $u$ 
and $v$, respectively.
\end{proof}

We now prove that $I_t(L)$ defines an arithmetically Cohen-Macaulay variety. 

\begin{thm}\label{cm}
Mixed ladder determinantal ideals are Cohen-Macaulay.
\end{thm}

\begin{proof}
Let $I_t(L)$ be a mixed ladder determinantal ideal. We follow the 
Notation~\ref{ladder}. Let $t_{\mmax}, t_{\mmin}, L_{\mmax}, L_{\mmin}$ be
as in Proposition~\ref{maxmin}. Let $I_{\mmax}=I_{t_{\mmax}}(L_{\mmax})$ and 
$I_{\mmin}=I_{t_{\mmin}}(L_{\mmin})$.
We saw that 
$$K[L_{\mmax}]/I_{\mmax}[u^{-1}]\cong 
K[L]/I_t(L)[w_1,\ldots,w_i][w_1^{-1},\ldots,w_{i_0}^{-1}],$$
and 
$$K[L_{\mmin}]/I_{\mmin}[z_1,\ldots,z_j][z_1^{-1},\ldots,z_{j_0}^{-1}]\cong 
K[L]/I_t(L)[v^{-1}].$$
The rings $K[L_{\mmax}]/I_{\mmax}$ and $K[L_{\mmin}]/I_{\mmin}$ are Cohen-Macaulay by
Corollary~4.10 in~\cite{he92a}. So the thesis follows from the observation that a ring $S$ 
is Cohen-Macaulay if and only if $S[x]$ is Cohen-Macaulay.
\end{proof}

We now wish to characterize the arithmetically Gorenstein mixed
ladder determinantal varieties.
In Theorem~5.2 of~\cite{co96a}, Conca classifies the 
arithmetically Gorenstein ladder determinantal varieties. 
Using the result of Conca, Gonciulea and Miller classify the 
arithmetically Gorenstein mixed ladder determinantal varieties
associated to one-sided ladders (see Theorem~6.3.1
of~\cite{go00}).
Using the same result of Conca together with Proposition~\ref{maxmin}, 
we can characterize the arithmetically Gorenstein mixed ladder
determinantal varieties. In analogy with the paper~\cite{co96a}, we start by
the following observations.

\begin{rmks}\label{props}
a) Let $L$ be a {\bf disconnected} ladder, i.e. $L=L_1\cup L_2$
with $L_1\cap L_2=\emptyset$. Let $(d_1,c_1),\ldots,(d_j,c_j)$ be
the lower outside corners of $L_1$ and let
$(d_{j+1},c_{j+1}),\ldots,(d_k,c_k)$ be the lower outside corners
of $L_2$. Let $\sigma=(t_1,\ldots,t_j)$ and
$\tau=(t_{j+1},\ldots,t_k)$.
Then $$K[L]/I_t(L)\cong
K[L_1]/I_{t^{(1)}}(L_1)\otimes_K K[L_2]/I_{t^{(2)}}(L_2).$$
Hence $K[L]/I_t(L)$ is Gorenstein if and only if both
$K[L_1]/I_{\sigma}(L_1)$ and $K[L_2]/I_{\tau}(L_2)$ are
Gorenstein.

b) Let $L$ be a ladder such that for some $i\in\{1,\ldots,h-1\}$
and some $j\in\{1,\ldots,k-1\}$ the following inequalities are
satisfied
\begin{equation}\label{ineqGor}
b_{i+1}\geq d_j-t_j+2,\;\;\;\; a_i\leq c_{j+1}+t_{j+1}-2.
\end{equation}
Let $X_1$ be the submatrix of $X$ consisting of the first $d_j$
rows and the first $a_i$ columns. Let $X_2$ be the submatrix of 
$X$ consisting of the rows $b_{i+1},\ldots,m$
and the columns $c_{j+1},\ldots,n$. If $X_1\cap
X_2\neq\emptyset$, we change the names of
the indeterminates in $X_2$ so that $X_1\cap
X_2=\emptyset$. E.g. we can let $$X_1=(x_{v,u})_{1\leq v\leq
  d_j;\; 1\leq u\leq a_i} \;\;\;\mbox{and}\;\;\;
X_2=(y_{v,u})_{b_{i+1}\leq v\leq m;\; c_{j+1}\leq u\leq n}.$$
Let $L_1$ be the subladder of $X_1$ with lower outside corners
$(d_1,c_1),\ldots,(d_j,c_j)$ and upper outside corners
$(b_1,a_1),\ldots,(b_i,a_i)$. Let $L_2$ be the subladder of $X_2$ with lower outside corners
$(d_{j+1},c_{j+1}),\ldots,(d_k,c_k)$ and upper outside corners
$(b_{i+1},a_{i+1}),\ldots,(b_h,a_h)$. Let $\sigma=(t_1,\ldots,t_j)$ and
$\tau=(t_{j+1},\ldots,t_k)$.
Denote by $R_t(L)$ the
quotient $K[L]/I_t(L)$, similarly for $L_1$ and $L_2$. Then
$$R_t(L)\cong R_{\sigma}(L_1)\otimes_K R_{\tau}(L_2)/I$$ where
$I$ is generated by $(d_j-b_{i+1}+1)(a_i-c_{j+1}+1)$ linear
forms. We have $$\hgt I=\dim [R_{\sigma}(L_1)\otimes_K
R_{\tau}(L_2)] - \dim R_t(L) = (d_j-b_{i+1}+1)(a_i-c_{j+1}+1)$$
hence $I$ is generated by a regular sequence. As in part a) we
conclude that $K[L]/I_t(L)$ is Gorenstein if and only if both
$K[L_1]/I_{\sigma}(L_1)$ and $K[L_2]/I_{\tau}(L_2)$ are
Gorenstein.

c) It follows from a) and b) that for the purpose of characterizing the
arithmetically Gorenstein mixed ladder determinantal varieties,
we can assume without loss of generality that $L$ is
connected. Moreover, we can assume that for each pair $i,j$ at least one of the 
inequalities (\ref{ineqGor}) is not satisfied. 
In fact, if $L$ is either disconnected or the
inequalities (\ref{ineqGor}) hold for some $i,j$, then we can write $L=L_1\cup L_2$ 
and reduce to testing the Gorensteinness on the subladders $L_1$ and $L_2$. 
Moreover, we can assume without loss of generality that 
$t_{\mmin}=\mmin\{t_2,\ldots,t_k\}\geq 2$.
\end{rmks}

We are now ready to give a numerical criterion to decide whether the ring 
$K[L]/I_t(L)$ is Gorenstein, for a ladder $L$. We assume without
loss of generality that $L$ satisfies the properties listed in part
c) of Remarks~\ref{props}.

\begin{thm}\label{ag}
Let $L$ be a ladder satisfying the properties in
Remarks~\ref{props}, part c).
Let 
$$J=\{j\;|\; c_j+t_j=c_{j+1}+t_{j+1}\}\subseteq\{1,\ldots,k-1\}$$
and
$$H=\{j+1\;|\;
d_j-t_j=d_{j+1}-t_{j+1}\}\subseteq\{1,\ldots,k-1\}.$$
Let $\{u_1,\ldots,u_s\}=\{1,\ldots,k\}\setminus(H\cup J).$
Then $I_t(L)$ is Gorenstein if and only if all of the following hold:
\begin{itemize}
\item $m-t_k=n-t_1$,
\item $c_{u_{j+1}}-d_{u_j}=2+t_{u_1}-t_{u_j}-t_{u_{j+1}}\;\; \mbox{for}\;\;
  j\in\{1,\ldots,s-1\}$
\item $a_i-b_{i+1}=t_1-2\;\;\mbox{for}\;\; i=1,\ldots,h-1$.
\end{itemize}
\end{thm}

\begin{proof}
Let $I_t(L)$ be a mixed ladder determinantal ideal. We follow the 
Notation~\ref{ladder}. Let $L_{\mmin},L_{\mmax},I_{\mmin},I_{\mmax}$
be as in Proposition~\ref{maxmin}. We assume without loss of
generality that $t_{\mmin}\geq 2$ (see item c) of
Remarks~\ref{props}).
To prove necessity of the conditions, assume that $K[L]/I_t(L)$ is a
Gorenstein ring.
By Proposition~\ref{maxmin} we have 
$$K[L_{\mmin}]/I_{\mmin}[z_1,\ldots,z_j][z_1^{-1},\ldots,z_{j_0}^{-1}]\cong 
K[L]/I_t(L)[v^{-1}]$$
therefore $K[L_{\mmin}]/I_{\mmin}$ is Gorenstein.
By Theorem~5.2 of~\cite{co96a}, the ring $K[L_{\mmin}]/I_{\mmin}$ 
is Gorenstein if and only if 
the smallest submatrix containing $L_{\mmin}$ is a square matrix, 
and the inside lower corners 
(resp. inside upper corners) lie on 
certain diagonals. The assumption that we need in order to apply
the theorem to the ladder $L_{\mmin}$ are satisfied, since they can be
rewritten as the assumptions in part b) of Remarks~\ref{props}. 
The fact that at least one of the inequalities
(\ref{ineqGor}) is not satisfied corresponds to the assumption that 
no outside lower corner of $L_{\mmin}$ belongs to the upper border of
$L_{\mmin}$ (i.e., no outside lower corner of $L_{\mmin}$ 
belongs to the same row or column of an outside upper corner).
The smallest submatrix containing $L_{\mmin}$ has $m-t_k+t_{\mmin}$ rows and 
$n-t_1+t_{\mmin}$ columns. Therefore it is square if and only if $m-t_k=n-t_1$.
The numerical conditions that characterize Gorensteinness in
\cite{co96a} translate into the last two conditions 
in the thesis. 
In particular, we can eliminate the lower outside corners of
$L_{\mmin}$ corresponding to the indexes in $H\cup J$. In fact, if two 
corners of $L_{\mmin}$ belong to the same row, we can eliminate the 
rightmost one. If two corners of $L_{\mmin}$ 
belong to the same column, then we can eliminate the one in the row above.
Therefore, we need
to check condition (2) only on a subset of the lower outside corners.
So the conditions are necessary for the ideal $I_t(L)$ to be
Gorenstein.

To show sufficiency of the conditions, assume that they are satisfied.
By Proposition~\ref{maxmin} we have 
$$K[L_{\mmax}]/I_{\mmax}[u^{-1}]\cong K[L]/I_t(L)[w_1,\ldots,w_i]
[w_1^{-1},\ldots,w_{i_0}^{-1}].$$
Again, the conditions in the statement of the theorem 
translate into the numerical conditions that characterize 
Gorensteinness in~\cite{co96a}. Therefore
$K[L_{\mmax}]/I_{\mmax}[u^{-1}]$ is Gorenstein.
Moreover, the smallest submatrix containing $L_{\mmax}$ has
$m-t_k+t_{\mmax}$ rows and 
$n-t_1+t_{\mmax}$ columns. Therefore it is square, since by assumption 
$m-t_k=n-t_1$.
From the isomorphism above we conclude that $K[L]/I_t(L)$ is Gorenstein.
\end{proof}

\section{Linkage of mixed ladder determinantal varieties}

In this section we prove that any mixed ladder determinantal variety can be  
obtained from a linear variety by a finite sequence of ascending elementary 
G-biliaisons. The biliaisons are performed on mixed ladder determinantal 
varieties, which are irreducible, hence generically Gorenstein. 
Therefore we conclude that mixed ladder determinantal varieties are glicci.
In the end of the section, we sketch how one can show that the same conclusions 
hold for a mixed ladder determinantal variety associated to a matrix of generic 
forms. 

The next theorem is the main result of this paper. It is analogous 
to Gaeta's Theorem, to Theorem~3.6 of
\cite{kl01}, and to Theorem~4.1 of \cite{ha04u2}. It is also analogous to 
Theorem~2.3 of~\cite{go05u} for a matrix that is not symmetric.
The idea of the proof is as follows: starting from a mixed ladder determinantal 
variety $V$, we construct two more mixed ladder determinantal varieties $V'$ and $W$ 
such that $V$ and $V'$ are generalized divisors on $W$. We show that $V'$ can be 
obtained from $V$ by an elementary G-biliaison on $W$.

\begin{thm}\label{bilmix}
Any mixed ladder determinantal variety in $\PP^r$ can be obtained from a linear
variety by a finite sequence of ascending elementary biliaisons.
\end{thm}

\begin{proof}
Let $V\subset\PP^r=Proj(K[L])$ be a mixed ladder determinantal variety. We
follow the Notation~\ref{ladder}. Let $L\subseteq X=(x_{i,j})$ be
the ladder whose minors of size $t=(t_1,t_2,\ldots,t_k)$ define 
the variety $V$. In other words, 
$I_V=I_{t_1}(L_1)+\cdots+I_{t_k}(L_k)\subseteq R=K[L]$.

If $t_1=\ldots=t_k=1$ then $V$ is a linear variety. Therefore, we concentrate 
on the case when $t_i\geq 2$ for some $i$.

Let $L'$ be the subladder of $L$ defined by the lower outside corners
$$\{(d_j-t_j+1,c_j+t_j-1)\;|\; j=1,\ldots,k\}.$$
It follows from Theorem~\ref{codim}
that $I_t(L)$ defines a variety of codimension $c$ equal to the number 
of entries in $L'$.
Fix an $i$ such that $t_i\geq 2$ and let $(d_i,c_i)$ be the position of 
the $i$-th lower outside corner of $L$. 
Let $$t'=(t_1,\ldots,t_{i-1},t_i-1,t_{i+1},\ldots,t_k),$$ and
let $M$ be the subladder of $L$ defined by the lower outside corners 
$\{(\delta_j,\gamma_j)\;|\; j=1,\ldots,k\}$:
$$(\delta_j,\gamma_j)=\left\{\begin{array}{ll} 
(d_i-1,c_i+1) & \mbox{if}\;\; j=i, \\
(d_j,c_j) & \mbox{if}\;\; j\neq i.
\end{array}\right.$$ 
Informally, $M$ is obtained from $L$ by removing 
the last row and the first column from the $i$-th lower step 
of the ladder $L$ (in the case of consecutive lower outside corners that lie on 
the same row or column however, one has to be careful about the corners while canceling).
Notice that $M$ is a ladder according to Definition~\ref{ladder}, unless either 
$d_i=d_{i-1}$ or $c_i=c_{i+1}$. By Assumption~\ref{assumpt}, if $d_i=d_{i-1}$
then $t_{i-1}>t_i\geq 2$. Therefore we can replace $i$ by $i-1$ and construct a new 
subset $M$ of $L$, starting from the $(i-1)$-th lower outside corner of $L$. The process 
eventually ends, since the subset $M$ constructed started from the leftmost lower outside 
corner of $L$ on the row $d_i$ is always a ladder. If $c_i=c_{i+1}$, then $t_{i+1}>t_i\geq 2$. 
In this case  we can replace $i$ by $i+1$ and construct a new 
subset $M$ of $L$, starting from the $(i+1)$-th lower outside corner of $L$. The process 
eventually ends, since the subset $M$ constructed started from the lowest outside lower
corner of $L$ on the column $c_j$ is always a ladder. 
Summarizing, we have shown that starting from a ladder 
$L$ with $t_j\geq 2$ for some $j$, we can construct a ladder $M\subset L$, with 
$t'_i=t_i-1$ for some $i$ and $t'_j=t_j$ for $j\neq i$.
 
Again, by Theorem~\ref{codim} $I_{t'}(M)\subset K[M]$ defines a mixed 
ladder determinantal variety of codimension $c=\hgt I_1(L')$.
Hence $I_{t'}(M)\subset K[L]$ defines a variety $V'\subset\PP^r$ of codimension $c$,
which is a cone over a mixed ladder determinantal variety. From now on, we will not 
distinguish between mixed ladder determinantal varieties and cones over them.

Let $N$ be the ladder obtained from $L$ by removing 
the entry $x_{d_i,c_i}$. $N$ is a ladder, since both $L$ and $M$ are.
$N$ has the same upper outside corners as $L, M$, and it has lower outside 
corners $$\{(d_1,c_1),\ldots,(d_{i-1},c_{i-1}),(d_i-1,c_i),(d_i,c_i+1),
(d_{i+1},c_{i-1}),\ldots,(d_k,c_k)\}.$$
Let $$\tau=(t_1,\ldots,t_{i-1},t_i,t_i,t_{i+1},\ldots,t_k)$$
be the vector obtained from $t$ by repeating twice the entry $t_i$.
Let $W$ be the mixed ladder 
determinantal variety with associated ideal $I_{\tau}(N)$. 
It follows from Theorems~\ref{codim}, ~\ref{prime}, and ~\ref{cm}
that $W$ has codimension $c-1$, and is irreducible and arithmetically 
Cohen-Macaulay. So $W$ is generically Gorenstein and $V, V'$ are 
generalized divisor on it. 

We denote by $H$ a hyperplane section divisor on $W$.
We claim that $$V\sim V'+H,$$ 
where $\sim$ denotes
linear equivalence of divisors on $W$. It follows that $V$ is
obtained by an ascending elementary biliaison from $V'$.
Repeating this argument, after $t_1+\cdots+t_k-k$ biliaisons we
reduce to the case $t_1=\ldots=t_k=1$.
Then the resulting variety is linear, 
hence we have bilinked $V$ to a linear variety. Notice that the
ideal of this linear variety is generated by the entries of $L'$, where $L'$
is the subladder of $L$ defined in Theorem~\ref{codim}.

Let $\I_{V|W}$, $\I_{V'|W}$  be the ideal sheafs on $W$ of $V$ and
$V'$. In order to prove the claim we must show that 
\begin{equation}\label{isosh} 
\I_{V|W}\cong \I_{V'|W}(-1).
\end{equation}

A minimal system of generators of $I_{V|W}=H^0_*(\I_{V|W})=I_t(L)/I_{\tau}(N)$ is
given by the images in the coordinate ring of $W$ of the $t_i\times t_i$ minors 
of $L$ that involve the entry $x_{d_i,c_i}$.
To keep the notation simple, we denote both an element of $R$ and
its image in $R/I_W$ with the same symbol, and we set $t=t_i$, $d=d_i$, $c=c_i$.
We let a minor $L_{i_1,\ldots,i_{t-1},d;c,j_1,\ldots,j_{t-1}}$, 
resp. $L_{i_1,\ldots,i_{t-1};j_1,\ldots,j_{t-1}}$, be $0$ if it involves one 
or more entries that do not belong to $L$. Hence we have
$$I_{V|W}=(L_{i_1,\ldots,i_{t-1},d;c,j_1,\ldots,j_{t-1}} \; | \; 1\leq
i_1<\ldots<i_{t-1}<d, c<j_1<\ldots<j_{t-1}\leq n).$$

A minimal system of generators of $I_{V'|W}=H^0_*(\I_{V'|W})=
I_{t-1}(M)/I_t(N)$ is
given by the images in the coordinate ring of $W$ of the minors
of $L_i$ of size $t_i-1$ that do not involve row $d$ and column $c$:
$$I_{V'|W}=(L_{i_1,\ldots,i_{t-1};j_1,\ldots,j_{t-1}} \; |\;
1\leq i_1<\ldots<i_{t-1}<d, c<j_1<\ldots<j_{t-1}\leq n).$$
Notice that with our conventions we have 
$$L_{i_1,\ldots,i_{t-1},d;c,j_1,\ldots,j_{t-1}}=0 \;\;\; \mbox{if and only if}
\;\;\;L_{i_1,\ldots,i_{t-1};j_1,\ldots,j_{t-1}}=0$$
where $L_{i_1,\ldots,i_{t-1},d;c,j_1,\ldots,j_{t-1}}$ is a generator of $I_{V|W}$
and $L_{i_1,\ldots,i_{t-1};j_1,\ldots,j_{t-1}}$ is a generator of $I_{V'|W}$.

In order to produce an isomorphism as in (\ref{isosh}), it suffices to
observe that the ratios 
\begin{equation}\label{eqk}
\frac{L_{i_1,\ldots,i_{t-1},d;c,j_1,\ldots,j_{t-1}}}
{L_{i_1,\ldots,i_{t-1};j_1,\ldots,j_{t-1}}}
\end{equation}
are all equal as elements of $H^0(\K_W(1))$, where $\K_W$ is the
sheaf of total quotient rings of $W$. Then the isomorphism
(\ref{isosh}) is given by multiplication by that element. 

Equality of all the ratios in (\ref{eqk}) follows from Lemmas~2.4 and ~2.6
in~\cite{go05u}, where we prove that
$$L_{i_1,\ldots,i_{t-1},d;c,j_1,\ldots,j_{t-1}}\cdot
L_{k_1,\ldots,k_{t-1};l_1,\ldots,l_{t-1}}-
L_{k_1,\ldots,k_{t-1},d;c,l_1,\ldots,l_{t-1}}\cdot
L_{i_1,\ldots,i_{t-1};j_1,\ldots,j_{t-1}}\in
I_W$$
for any choice of $i_1,\ldots,i_{t-1}$, $j_1,\ldots,j_{t-1}$, 
$k_1,\ldots,k_{t-1}$, $l_1,\ldots,l_{t-1}$. This completes the proof of the 
claim and of the theorem. 
\end{proof}

Standard results in liaison theory give the following.

\begin{cor}\label{glicci}
Every mixed ladder determinantal variety $V$ can be bilinked in $t_1+\ldots+t_k-k$ 
steps to a linear variety of the same codimension.
In particular determinantal varieties associated to matrices of indeterminates, 
ladder determinantal varieties, and mixed ladder determinantal varieties are glicci.
\end{cor}

\begin{rmk}
In Theorem~\ref{ag} we classify the mixed ladder determinantal varieties which are 
arithmetically Gorenstein. The fact 
that arithmetically Gorenstein schemes are glicci follows 
from a result of Casanellas, Drozd, and~Hartshorne (see Theorem~7.1 of~\cite{ca04u}). 
\end{rmk}

We conclude by noticing that by the upper-semicontinuity principle, 
one can argue that a generic 
mixed ladder determinantal scheme is glicci. The idea is as follows:
fix $m,n$ and upper and lower outside corners for a ladder $L$, and 
assume that all the ideals involved in the proof of Theorem~\ref{bilmix}
have maximal height (i.e. the height of the ideal that we would obtain 
if all the entries in $L$ were distinct indeterminates). This is the case
e.g. when the entries of $L$ are generic polynomials in $S=K[x_0,\ldots,x_s]$
and $s\gg 0$. Let $T$ be the scheme with saturated ideal 
$I_T=I_t(L)\subseteq S$. We call such a scheme a 
{\bf generic mixed ladder determinantal scheme}. 
The proof of Theorem~\ref{bilmix} can be repeated with no changes for $T$, 
under the assumption that $T$ is a generic mixed ladder determinantal 
scheme. Therefore we have the following.

\begin{thm}
Any generic mixed ladder determinantal scheme can be obtained from a 
complete intersection by a finite sequence of ascending elementary G-biliaisons.
In particular, a generic mixed ladder determinantal scheme is glicci.
\end{thm}

\end{document}